\documentclass[11pt]{article}   %% Standard LaTeX.
\usepackage{amsmath,amssymb}    %% For better support of math.
\usepackage{pstricks}                                
\usepackage{amsthm}                               
\usepackage{url}                %% Supports formating URLs.
\usepackage{graphicx}           %% For eps figures
\usepackage{subfigure}

\begin{document}

\title{AN IMPROVED LOWER BOUND FOR MOSER'S WORM PROBLEM}

\author{TIRASAN KHANDHAWIT AND SIRA SRISWASDI}
\date{}     
\maketitle

\theoremstyle{plain}% default
\newtheorem{thm}{Theorem}
\newtheorem{cor}[thm]{Corollary}
\newtheorem{prop}[thm]{Proposition}
\newtheorem{lem}[thm]{Lemma}
\newtheorem{conj}{Conjecture}
\newtheorem{quest}{Question}

\theoremstyle{definition}
\newtheorem{exam}[thm]{Example}
\newtheorem{defin}[thm]{Definition}

\theoremstyle{remark}
\newtheorem*{rem}{Remark}
\newtheorem*{note}{Note}
\newtheorem*{proof}{Proof}

\begin{abstract}
We show that any convex region which contains a unit segment, an equilateral triangle of sides $\frac{1}{2}$, and a square of side $\frac{1}{3}$ always has area at least 0.227498. Using grid-search algorithm, we attempt to find a configuration of these three objects with minimal convex hull area. Consequently, we improve a lower bound for Moser's worm problem from 0.2194 to 0.227498.
\end{abstract}

\section{Introduction}\label{sec-intro}	%% Refer to via \ref{sec-intro}.
In 1966, Leo Moser \cite{Moser} asked for the region of smallest area which can accommodate every planar arc of length one. The problem is known as ``Moser's worm problem'' and is a variation of universal cover problems (see \cite{rpdg}). In Moser's problem, a cover is a set which contains a copy of any rectifiable planar arc of unit length, and is usually assumed to be convex. Such a minimal cover is known to have area between 0.2194 and 0.2738. However, the original problem remains unsolved.  

There have been many works to find a universal cover for any unit planar arc. The first few such covers are Meir's semidisc of diameter 1 with an area of 0.39269 in and John Wetzel's sectorial plate with an area of 0.34510 \cite{Wet73}. In 1974, Gerriets and Poole \cite{jg} introduced a rhombus cover with an area of 0.28870. In 2003, Norwood and Poole \cite{rg} constructed a non-convex cover of area 0.260437 whose convex hull gave the current best upper bound of 0.2738. Furthermore, Wetzel \cite{rgm} has conjectured an upper bound of 0.23450.

On the other hand, there have not been much improvement for a lower bound. Wetzel \cite{Wet73} gave the lower bound of 0.2194 in 1973 by using Schaer's broadworm \cite{Sch}, a unit arc whose width is at least 0.4389 in every direction. 

To improve a lower bound, we observe that any convex cover must contain a unit segment, an equilateral triangle of sides $\frac{1}{2}$, and a square of side $\frac{1}{3}$. We then study all possible configurations of these three objects, i.e. a placement of the three objects in the plane. Our main result is stated below.

\begin{thm} \label{main} For any configuration of a unit segment, an equilateral triangle of sides $\frac{1}{2}$, and a square of side $\frac{1}{3}$, its convex hull always has area at least 0.227498.
\end{thm}

It follows immediately that.

\begin{thm} A convex universal cover for a unit length arc has area at least 0.227498.
\end{thm}

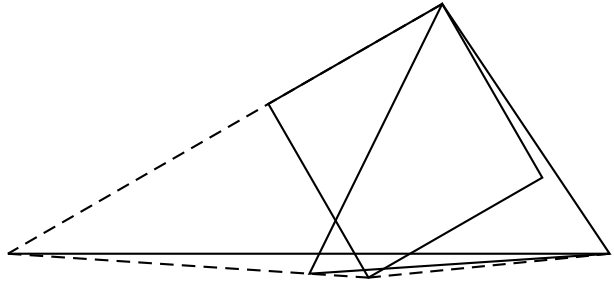
\begin{figure}[htbp]
\centering
\begin{pspicture}(-4,-.4)(4,3.6) \psset{xunit=0.8cm,yunit=0.8cm, runit=0.8}
  \psline(-5,0)(5,0)
  \pspolygon(2.2177,4.1543)(0.0111,-0.3324)(5,0)
  \pspolygon(2.2177,4.1543)(-0.6713,2.4915)(0.9915,-0.3975)(3.8805,1.2654)
  \psline[linestyle=dashed](2.2177,4.1543)(-5,0)
  \psline[linestyle=dashed](-5,0)(0.9915,-0.3975)
  \psline[linestyle=dashed](5,0)(0.9915,-0.3975)
\end{pspicture}
  \caption{Current best configuration with area of 0.227589669377.}\label{champ}
\end{figure}

We note that Brass and Sharifi \cite{leblower} made a similar observation to improve a lower bound for Lebesgue's universal cover problem. In this paper, we use geometric and some analytic argument to prove the main theorem in section~\ref{pf}. In section~\ref{search}, we outline the heuristic grid-search algorithm similar to \cite{leblower} to search for small configurations and find a configuration with an area of 0.227589669377 (Figure~\ref{champ}).

\textbf{Acknowledgments.} The authors are indebted to W. Wichiramala for his helpful advises and supervision. We would also like to thank colleagues at Chulalongkorn University for useful discussions, especially R. Tanadkithirun who found a small configuration closed to the current best configurattion.

\section{Proof of the Main Theorem}\label{pf}
Let $\mathcal{U}$ be a set containing all configuration of a unit segment $\mathcal{L}$, an equilateral triangle $\mathcal{T}$ of side $\frac{1}{2}$, and a square $\mathcal{S}$ of side $\frac{1}{3}$ in a Euclidean plane. With out loss of generality, we can assume that $\mathcal{L}$ is a segment from $(0,0)$ to $(1,0)$.

A square $\mathcal{S}$ can be described by three parameters $(x_1,y_1,\alpha)$ where $(x_1,y_1)$ is a center of $\mathcal{S}$ and $\alpha$ is an angle of rotation. This means that vertices of $\mathcal{S}$ are $(x_1 + \frac{\sqrt{2}}{6} \cos{\alpha},y_1 + \frac{\sqrt{2}}{6} \sin{\alpha})$ , $(x_1 + \frac{\sqrt{2}}{6} \cos{(\alpha + \pi/2)},y_1 + \frac{\sqrt{2}}{6} \sin{(\alpha + \pi/2)})$, $(x_1 + \frac{\sqrt{2}}{6} \cos{(\alpha + \pi)},y_1 + \frac{\sqrt{2}}{6} \sin{(\alpha + \pi)})$, and \\ $(x_1 + \frac{\sqrt{2}}{6} \cos{(\alpha + 3\pi/2)},y_1 + \frac{\sqrt{2}}{6} \sin{(\alpha + \pi/2)})$. By rotational symmetry of a square, we can assume $0 \leq \alpha \leq \pi/2$. Similarly we describe $\mathcal{T}$ by parameters $(x_2,y_2,\beta)$ with $0 \leq \beta \leq 2\pi/3$.

Notice that a reflection across a line $x = \frac{1}{2}$ sends a square given by parameters $(x_1,y_1,\alpha)$ to $(\frac{1}{2} - x_1,y_1,\pi/2 - \alpha)$. We can then assume without loss of generality that $\pi/4 \leq \alpha \leq \pi/2$. We also notice that a half turn centered at $(\frac{1}{2},0)$ fixes a square while sending a triangle $(x_2,y_2,\beta)$ to $(1 - x_2,-y_2,\beta + \pi/3)$. Hence we can assume further that $\pi/3 \leq \beta \leq 2\pi/3$.

We now define a map $\phi : \mathbb{R}^6 \rightarrow \mathcal{U}$ by sending $(x_1,y_1,\alpha,x_2,y_2,\beta)$ to a square and a triangle with parameters $(x_1,y_1,\alpha)$ and $(x_2,y_2,\beta)$ respectively. Note that $\phi$ is surjective but not injective. A configuration determines centers of a square and a triangle uniquely but the angles are determined up to rotational symmetry. For each configuration $X$ of the three objects, let $\mathcal{C}(X)$ denote its convex hull and $\mu(X)$ denote the area of $\mathcal{C}(X)$. It is clear that the composition $\mu \circ \phi : \mathbb{R}^6 \rightarrow \mathbb{R}$ is continuous.

Before estimating area, we introduce a notion of height.

\begin{defin} Let $\vec{u}$ and $\vec{v}$ be vectors in $\mathbb{R}^2$. The \emph{height} of $\vec{v}$ with respect to $\vec{u}$ is a value $\left|\vec{v}\right| \sin{\theta}$ denoted by $h_{\vec{u}}(\vec{v})$, where $\theta$ is an angle between $\vec{u}$ and $\vec{v}$.   
\end{defin}

Label points $(0,0)$ and $(1,0)$ with E and F respectively. Label vertices of a square $\mathcal{S}$ with A, B, C, D starting from $(x_1 + \frac{\sqrt{2}}{6} \cos{\alpha},y_1 + \frac{\sqrt{2}}{6} \sin{\alpha})$ counterclockwise. Label vertices of a triangle $\mathcal{T}$ with P,Q,R starting from $(x_2 + \frac{\sqrt{3}}{6} \cos{\beta},y_2 + \frac{\sqrt{3}}{6} \sin{\beta})$ counterclockwise. We have the following inequalities.

\begin{lem} \label{1stineq} Let $X = \phi(x_1,y_1,\alpha,x_2,y_2,\beta)$ be a configuration. Then

(1) $\mu(X) \geq \frac{\sqrt{2}}{6} \sin{\alpha}$

(2) $\mu(X) \geq \mbox{Max} \left\{ \frac{1}{4} \sin{(\beta - \pi/6)}, \frac{1}{4} \sin{(\beta + \pi/6)} \right\}$
\end{lem}  

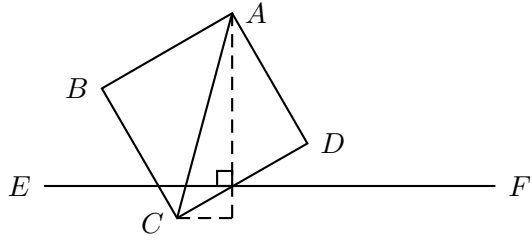
\begin{figure}[htbp] 
\centering \begin{pspicture}(-3,-.5)(3,2.5)
  \psline(-3,0)(3,0)
  \pspolygon(0.5,0.566022)(-0.5,2.2980762)(-2.232052,1.29807)(-1.232052,-0.426716)
  \uput[l](-3,0){$E$}
  \uput[r](3,0){$F$}
  \uput[l](-1.232052,-0.5){$C$}
  \uput[r](0.5,0.566022){$D$}
  \uput[r](-0.5,2.2980762){$A$}
  \uput[l](-2.232052,1.29807){$B$}
  \psline(-1.232052,-0.426716)(-0.5,2.2980762)
  \psline[linestyle=dashed](-.5,-0.426716)(-0.5,2.2980762) \psline[linestyle=dashed](-.5,-0.426716)(-1.232052,-0.426716)
  \pspolygon(-.5,0)(-.7,0)(-.7,.2)(-.5,.2)
  \end{pspicture}
  \caption{A square $\cal{S}$ and a height of AC with respect to EF.} \label{sqheight}
\end{figure}

\begin{proof} (1) Consider a vector $\overrightarrow{CA}$ , which is described by $(\frac{\sqrt{2}}{3} \cos{\alpha}, \frac{\sqrt{2}}{3} \sin{\alpha})$. Since $\pi/4 \leq \alpha \leq \pi/2$, the height $h_{\overrightarrow{EF}}(\overrightarrow{CA})$ is $\frac{\sqrt{2}}{3} \sin{\alpha}$. Then $\mu(X) \geq \mu(AECF) \geq \frac{1}{2} \left|\overrightarrow{EF}\right| h_{\overrightarrow{EF}}(\overrightarrow{CA}) = \frac{\sqrt{2}}{6} \sin{\alpha}$ (See figure~\ref{sqheight}).

(2) Similarly, $\mu(X) \geq \mu(PEQF) \geq \frac{1}{4} \sin{(\beta - \pi/6)}$ and \\ $\mu(X) \geq \mu(PFRE) \geq \frac{1}{4} \sin{(\beta + \pi/6)}$ as $\pi/3 \leq \beta \leq 2\pi/3$.

\end{proof}

From this point, we will try to reduce the domain $\mathbb{R}^6$ to a smaller subset on which it is still sufficient for us to search for the smallest cover. We will first try to exclude a configuration of which convex area is greater than 0.23. Let $\mathcal{K'}_1$ be a subset $\mathbb{R}^2 \times \left[45^\circ,78^\circ \right] \times \mathbb{R}^2 \times \left[83^\circ,97^\circ \right] $ of $\mathbb{R}^6$. We compute that $\arcsin{(\frac{6}{\sqrt{2}} 0.23)} = 1.35 \approx 77.3^\circ $ and $\arcsin{(4 \times 0.23)} = 1.168 \approx 66.9^\circ$. By lemma~\ref{1stineq}, it follows that $\mu(X) > 0.23$ if $X \notin \phi(\mathcal{K}'_1)$. This gives a subdomain $\phi(\mathcal{K'}_1) \subset \mathcal{U}$ denoted by $\mathcal{K}_1$.
Next, we show that a global minimum of $\mu \circ \phi$ exists.

\begin{lem} A function $\mu \circ \phi : \mathbb{R}^6 \rightarrow \mathbb{R} $ attains its minimum value.
\end{lem}

\begin{proof}
A square $\mathcal{S}$ contains the inscribed circle of radius $\frac{1}{6}$. The convex hull of this circle and the origin has area at least $\sqrt{\frac{x_1^2 + y_1^2}{36} - \frac{1}{6^4}}$ , when $x_1^2 + y_1^2$ is greater than $ \frac{1}{36}$. We compute that if $x_1^2 + y_1^2$ is greater than 1.4, the area of is greater than $0.23$. Analogously, the convex hull of the inscribed circle of the triangle $\mathcal{T}$ and the origin is greater than $0.23$ when $x_2^2 + y_2^2$ is larger than 1.7. Since a subset $B_{1.4} \times \left[0,\pi/2\right] \times B_{1.7} \times \left[0,2\pi/3\right] $ , where $B_r$ denotes a disk centered at the origin in $\mathbb{R}^2$, is compact, the minimal value of $\mu \circ \phi$ is attained.
\qed \end{proof}    

We consider a subset $\mathcal{K}_2$ of $\mathcal{U}$ with the following properties: For a configuration in $\mathcal{K}_2$,

 (i) A distance between any point in $\mathcal{S}$ or $\mathcal{T}$ and $\mathcal{L}$ is not more than $1$.  

(ii) $\mathcal{S}$ and $\mathcal{T}$ lies in a region $-0.46 \leq y \leq 0.46$ 

(iii) Both $\mathcal{S}$ and $\mathcal{T}$ have non-empty intersection with $\mathcal{L}$ 

\begin{prop} If a configuration $X$ is not in $\mathcal{K}_2$, then $\mu(X)$ is not minimal.
\end{prop}

\begin{proof} (i) Suppose, without loss of generality, that there exist a point $(x_0,y_0)$ in $\mathcal{S} \cup \mathcal{T}$ with distance more than 1 from the origin. We can form a new configuration $X'$ consisting of $\mathcal{S}$, $\mathcal{T}$, and a segment from the origin to $\frac{1}{\sqrt{x_0^2 + y_0^2}}(x_0,y_0)$. We easily see that $\mathcal{C}(X)$ contains $\mathcal{C}(X')$ as a subset, and hence is not minimal.

(ii) Suppose that $\mathcal{S} \cup \mathcal{T}$ contains a point $(x_0,y_0)$ with $\left|y_0\right| > 0.46$. The area of a triangle with vertices $(0,0)$, $(1,0)$, and $(x_0,y_0)$ is $\frac{1}{2}\left|y_0\right| > 0.23$. Since there exists a configuration of which convex hull has area less than 0.23, $\mu(X)$ is not minimal.

(iii) Suppose that the square $\mathcal{S}$ lies above $\mathcal{L}$. Here we can further assume that $d(\mathcal{S} \cup \mathcal{T} , \mathcal{L}) \leq 1$, otherwise $\mu(X)$ is not minimal from above. This implies that $\mathcal{S}$ lies in a region $0 \leq x \leq 1$. Let $(x_0,y_0)$ be a point in $\mathcal{S}$ with minimal y-coordinate. We see that a translation $\mathcal{S} - (0,y_0)$ of $\mathcal{S}$ down by $y_0$ lies in $\mathcal{X}$. Thus $\mathcal{C}(X)$ contains $\mathcal{C}(X')$ as a subset, where $X'$ is a new configuration formed by $\mathcal{T}$, $\mathcal{L}$, and $\mathcal{S} - (0,y_0)$. The same argument applies to other cases. \qed     
\end{proof}

In order to prove the main theorem, we need another inequality. 

\begin{prop} Let $X$ be a configuration in $\mathcal{K}_1 \cap \mathcal{K}_2$ described by parameters $(x_1,y_1,\alpha,x_2,y_2,\beta)$. Then $\mu(X) \geq \frac{1}{6} (\frac{1}{2}\cos{(\alpha - \beta + 15^\circ)} + \cos{(\alpha - 45^\circ)}). $
\label{2ndineq} \end{prop}

\begin{proof} Let $\cal{W}$ be the strip of width $\frac{1}{3}$ bounded by extended segments $AB$ and $CD$ of $\cal{S}$ and $\cal{V}$ be the strip of bounded by extended segments $BC$ and $AD$ (See figure~\ref{fig-R}). We will consider several cases depending on the position of $E,F$ relative to $\cal{W}$, which has non-negative slope. The main argument is to estimate the area of a part of $\mathcal{C}(X)$ lying outside $\cal{S}$ by using sides of $\cal{S}$, the segment $EF$, and the side $PR$ of $\cal{T}$.

\begin{figure}[htbp]
\begin{center}
  \begin{pspicture}(-1.6,-1.5)(1.6,1.5) \psset{xunit=0.6cm,yunit=0.6cm, runit=0.2}
  \pspolygon(-0.84,-1.13)(1.14,-0.86)(0.86,1.12)(-1.12,0.85)
  \uput[dl](-0.84,-1.13){$C$}
  \uput[dr](1.19,-0.89){$D$}
  \uput[ur](0.86,1.16){$A$}
  \uput[ul](-1.12,0.89){$B$}
  \uput[r](1.44,0.27){\large $\cal{W}$} \uput[u](-.27,1.44){\large $\cal{V}$}
  \psline[linestyle=dashed](-2.82,-1.4)(3.12,-0.59)
  \psline[linestyle=dashed](-3.1,0.58)(2.84,1.39)
   \psline[linestyle=dashed](1.4,-2.82)(0.59,3.12)
  \psline[linestyle=dashed](-.58,-3.1)(-1.39,2.89)
\end{pspicture}
\end{center}
\caption{The strips $\cal{W}$ and $\cal{V}$ .}\label{fig-R}
\end{figure}
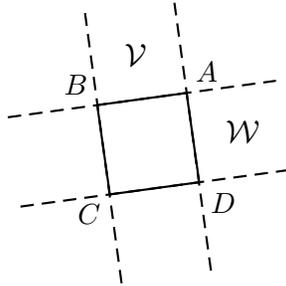

\underline{Case 1} Both E and F lie inside $\cal{W}$. We observe that a set $\mathcal{C}(X) - \mathcal{C}(\mathcal{S})$ contains either both triangles $BCE$ and $AFD$ or only one of the two triangles. This also depends on the position of $E$ and $F$ with respect to $\mathcal{V}$. Either both points lie outside $\mathcal{V}$ or one of the points lies outside $\mathcal{V}$. In the former case, a sum of height of two triangles with respect to $BC$ is exactly $\cos{(\alpha - 45^\circ)} - 1/3$. In the latter case, the triangle has height at least $\cos{(\alpha - 45^\circ)} - 1/3$ with respect to $BC$ (See figure~\ref{fig-case1}). Note that since $45^\circ \leq \alpha \leq 90^\circ$, we have $\cos{(\alpha - 45^\circ)} \geq 1/\sqrt{2} > 1/3$, the width of $\mathcal{V}$, and so at least one of $E$ and $F$ must lie outside $\mathcal{V}$.

\begin{figure}[htbp]
\begin{center}
  \begin{pspicture}(-3.1,-1.2)(3.1,1.8)
  \psline(-3,0)(3,0)
  \pspolygon(-0.84,-1.13)(1.14,-0.86)(0.86,1.12)(-1.12,0.85)
  \uput[r](3,0){$F$}
  \uput[l](-3,0){$E$} 
  \uput[l](-0.84,-1.13){$C$}
  \uput[r](1.19,-0.89){$D$}
  \uput[r](0.86,1.16){$A$}
  \uput[l](-1.12,0.89){$B$}
  \psline(-3,0)(-0.84,-1.13)
  \psline(-3,0)(-1.12,0.85)
  \psline(1.14,-0.86)(3,0)
  \psline(0.86,1.12)(3,0)
  \psline[linestyle=dashed](0.86,1.12)(-0.2,1.4)
   \psline[linestyle=dashed](-0.2,1.4)(-1.12,0.85)
   \psline[linestyle=dashed](0.2,-1.4)(-0.84,-1.13)
    \psline[linestyle=dashed](1.14,-0.86)(0.2,-1.4)
   \psline[linestyle=dashed,dash=1pt 6pt](-2.82,-1.4)(3.12,-0.59)
   \psline[linestyle=dashed,dash=1pt 6pt](-3.1,0.58)(2.84,1.39)
     \uput[d](0.2,-1.4){$R$}
     \uput[u](-0.2,1.4){$P$}
\end{pspicture}
\end{center}
\caption{Both $E$ and $F$ lie inside $\cal{W}$}\label{fig-case1}
\end{figure}
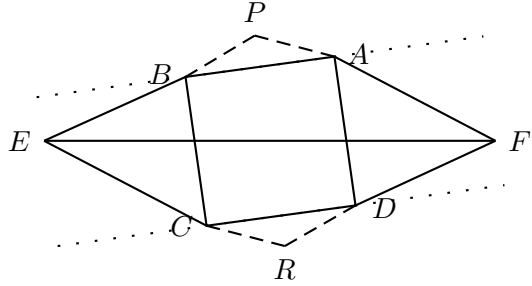

Similarly, a side $PR$ contributes to triangles above or below $\mathcal{S}$. Since $\alpha \in \left[45^\circ , 78^\circ\right]$ and $\beta \in \left[83^\circ , 97^\circ\right]$, we have $-37^\circ \leq \alpha - \beta + 15^\circ \leq 10^\circ$. We check that $\cos{37^\circ} = 0.798$ so we have $h_{AB}(PR) = \cos{(\alpha - \beta + 15^\circ)}/2 > 1/3$. Thus, at least one of P and R lies outside $\mathcal{W}$. We estimate the area of $\mathcal{C}(X)$ outside $\mathcal{S}$ by triangles of which bases are sides of the square (Figure~\ref{fig-case1}). Hence,

\begin{align*}
	\mu(X) & \geq  \frac{1}{9} + \frac{1}{6}(\cos{(\alpha - 45^\circ)} -  \frac{1}{3}) + \frac{1}{6}( \frac{1}{2}\cos{(\alpha - \beta + 15^\circ)} -  \frac{1}{3}) 
	\\ & \geq \frac{1}{6} (\frac{1}{2}\cos{(\alpha - \beta + 15^\circ)} + \cos{(\alpha - 45^\circ)})
\end{align*}

\underline{Case 2} E lies inside $\cal{W}$ and F lies below $\cal{W}$. We can suppose that $F$ lies on the right of $\mathcal{V}$, otherwise we can use the same argument as in Case 1. We notice that in this case the triangle $AFD$ may intersect with the triangle $CRD$ (See figure~\ref{fig-tri-ef}).

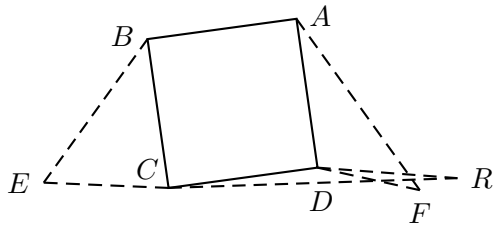
\begin{figure}[htbp]
\begin{center}
  \begin{pspicture}(-3.1,-1.2)(3.1,1.2)
  \pspolygon(-0.84,-1.13)(1.14,-0.86)(0.86,1.12)(-1.12,0.85)
  \uput[d](2.5,-1.16){$F$}
  \uput[l](-2.5,-1.06){$E$} 
  \uput[ul](-0.84,-1.13){$C$}
  \uput[d](1.19,-0.99){$D$}
  \uput[r](0.86,1.16){$A$}
  \uput[l](-1.12,0.89){$B$}
  \uput[r](3,-1){$R$}
  \psline[linestyle=dashed](-2.5,-1.06)(-0.84,-1.13)
  \psline[linestyle=dashed](-2.5,-1.06)(-1.12,0.85)
  \psline[linestyle=dashed](1.14,-0.86)(2.5,-1.16)
  \psline[linestyle=dashed](0.86,1.12)(2.5,-1.16)
    \psline[linestyle=dashed](1.14,-0.86)(3,-1)
  \psline[linestyle=dashed](-0.84,-1.13)(3,-1)
\end{pspicture}
\end{center}
\caption{Possible configuration in case 2.}\label{fig-tri-ef}
\end{figure}

To address the problem, we let $l_1$ be the line passing through $F$ parallel to $\cal{W}$ and $l_2$ be the line through $D$ and $F$ (see figure~\ref{fig-case2}). We can also assume that $R$ lies below $\mathcal{W}$, otherwise P lies above $\mathcal{W}$ and we can use the triangles $PAB$, $BEC$, and $AFD$ to estimate the area as the previous case. Consider two subcases:

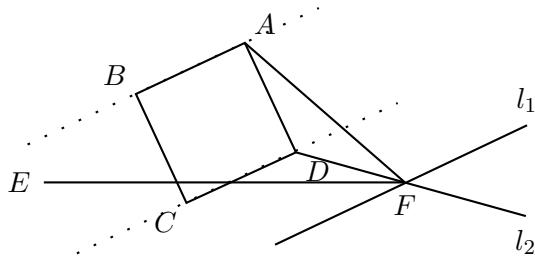
\begin{figure}
\centering
	\begin{pspicture}(-1.2,-.96)(4,2.2) \psset{xunit=0.8cm,yunit=0.8cm, runit=0.2}
  \psline(-3,0)(3,0)
  \pspolygon(-0.63,-0.34)(1.19,0.50)(0.34,2.32)(-1.47,1.47)
  \psline(5.01,-0.56)(1.19,0.50) %l2
  \psline(0.843,-1.036)(5.03,0.95) %l1
  \psline(0.34,2.32)(3,0) %CF
   \psline[linestyle=dashed,dash=1pt 6pt](-2.45,-1.18)(2.88,1.32)
   \psline[linestyle=dashed,dash=1pt 6pt](-3.27,0.63)(1.57,2.9)
  \uput[d](3,0){$F$}
  \uput[l](-3,0){$E$} 
  \uput[ur](0.34,2.32){$A$}
  \uput[ul](-1.47,1.47){$B$}
  \uput[dl](-.63,-.34){$C$}
  \uput[dr](1.19,.5){$D$}
  \uput[u](5.03,0.95){$l_1$}
  \uput[d](5.01,-0.56){$l_2$}
  \end{pspicture}
  \caption{The geometric construction in case 2}\label{fig-case2}
\end{figure}

\begin{itemize}
	\item $R$ lies above or on $l_1$. It follows that $h_{CD}(PR) \leq h_{CD}(PF)$, so we can use the height of $PF$ instead of $PR$ to estimate the area (i.e. use the triangle $CDF$ in place of $CDR$).
	\item $R$ lies below $l_1$. Since $X$ is in $\mathcal{K}_2$, the point $R$ must lie to the left of the line $x = 1$. Consequently, $R$ have to lie below $l_2$ as well. Then we see that the triangles $CRD$ and $AFD$ are disjoint. 
\end{itemize}

From both subcases, we still have the same inequality. 
	\[ \mu(X) \geq \frac{1}{6} (\frac{1}{2}\cos{(\alpha - \beta + 15^\circ)} + \cos{(\alpha - 45^\circ)}) \]

\underline{Case 3} E lies above $\cal{W}$ and F lies inside $\cal{W}$. Analogously, we can suppose that $E$ is on the left of $\cal{V}$ and $P$ lies above $\cal{W}$. Construct the line $l_1$ parallel to $AB$ at $E$ and the line $l_2$ joining $B$ and $E$. We now then use the same argument as in case 2 to obtain the result. 

\underline{Case 4} E lies above $\cal{W}$ and F lies below $\cal{W}$. In this case, we can directly apply the arguments from case 2 and case 3 together to prove the statement.

Finally, the cases that both E and F lie above or below $\cal{W}$ does not occur because $\cal{S}$ have to intersect $\cal{L}$ from a property of $\mathcal{K}_2$. 
\qed \end{proof}

Set
\begin{align*}
	f(\alpha,\beta) & = \frac{1}{6} (\frac{1}{2}\cos{(\alpha - \beta + 15^\circ)} + \cos{(\alpha - 45^\circ)}) \\
	g(\alpha) &= \frac{\sqrt{2}}{6} \sin{\alpha} \\
	h(\beta) &= \mbox{Max} \left\{ \frac{1}{4} \sin{(\beta - 30^\circ)}, \frac{1}{4} \sin{(\beta + 30^\circ)} \right\} \\
	p(\alpha,\beta) & = \mbox{Max} \left\{ f(\alpha,\beta), g(\alpha),h(\beta) \right\}
\end{align*}

We are now ready to prove the main theorem (Theorem~\ref{main}).

\begin{proof} From our definition, it suffices to consider an arbitrary configuration $X = \phi(x_1,y_1,\alpha,x_2,y_2,\beta)$ in $\mathcal{K}_1 \cap \mathcal{K}_2$. It follows from lemma~\ref{1stineq} and proposition~\ref{2ndineq} that $\mu(X) \geq p(\alpha,\beta)$. 

Suppose the contrary that there is $(\alpha_0,\beta_0)$ in $\left[45^\circ,78^\circ \right] \times \left[83^\circ,97^\circ \right] $ such that $p(\alpha_0,\beta_0) < 0.227498$. Using functions $g$ and $h$, we can further deduce that $\alpha_0 < 74.838^\circ$ and $84.496^\circ < \beta_0 < 95.504^\circ$, which implies $-35.504^\circ < \alpha_0 - \beta_0 + 15^\circ < 5.342^\circ$. As $\beta_0$ varies, we have $\cos{(\alpha_0 - \beta_0 + 15^\circ)} \geq \cos{(\alpha_0 - 80.504^\circ)}$. Then 
\begin{align*}
	f(\alpha_0,\beta_0) & \geq \frac{1}{6} (\frac{1}{2}\cos{(\alpha_0 - 80.504^\circ)} + \cos{(\alpha_0 - 45^\circ)}) \\
	& \geq 0.2274987 \ \ \mbox{as $\alpha_0 \in \left[45^\circ,74.838^\circ \right]$}
\end{align*}
which is a contradiction.
\qed \end{proof}

\section{Search for an Optimal Configuration} \label{search}
In this section, we search for a configuration with minimal area of convex hull to see how far its area is from our new lower bound.

The main strategy is to start off with large grid sizes $d_1$ for each $x_i$ and $y_i$ and $d_2$ for $\alpha$ and $\beta$. Then we heuristically zoom in, reducing the domain of candidate configurations, and reduce our grid sizes accordingly to gain more precision while keeping the computation time reasonable. First, we provide a theorem which relates the magnitude of grid sizes to the precision of our estimate for the area of the optimal configuration.

\begin{prop} \label{stepbound} Let $d_1$ and $d_2$ be grid sizes for xy-coordinate and angles in grid-search algorithm, respectively, then the error between the optimal area found in the search and the actual minimal area is at most $2.44916 \: d_1 + 0.49993 \: d_2$ .
\end{prop}

\begin{proof} Suppose $X = \phi(x_1,y_1,\alpha,x_2,y_2,\beta)$ is a configuration with minimal area of convex hull. Let $X' = \phi(x_1 + \delta_1,y_1 + \delta_2,\alpha + \theta_1,x_2 + \delta_3,y_2 + \delta_4,\beta + \theta_2)$ be a nearby configuration with $\left| \delta_i \right| \leq d_1 / 2$ and $\left| \theta_j \right| \leq d_2 / 2$ where $i = 1,\ldots,4$ and $j = 1,2$.   

We see that a vertex of the triangle of $X'$ can be obtained by a vertex of the triangle of $X$ via a rotation centered at $(x_1,y_1)$ with angle $\theta_1$ followed by translation in a direction $(\delta_1,\delta_2)$. Hence, the new vertex is far apart from the original vertex by distance at most $d_1 / \sqrt{2} + \sin{(d_2 /4)} / \sqrt{3} $. Similarly, a distance between vertices of the squares in $X$ and $X'$ is at most $d_1 / \sqrt{2} + \sqrt{2} \sin{(d_2 /4)} / 3$. Thus the vertices in $X'$ are at distance most $\delta = d_1 / \sqrt{2} + \sin{(d_2 /4)} / \sqrt{3} $ from vertices in $X$.

We now apply a Lipschitz bound used in \cite{leblower}.
\[ \left| \mu(X') - \mu(X) \right| \leq \delta \: \mbox{peri}(\mathcal{C}(X)) + \pi \delta^2 \]
Since $X$ lies in $\mathcal{K}_2$, the convex hull $\mathcal{C}(X)$ is contained in a convex region $\mathcal{D}$ bounded by lines $ y = \pm 0.46$ and circular arcs $x^2 + y^2 = 1$ and $(x-1)^2 + y^2 = 1$. We calculate that the perimeter of $\mathcal{D}$ is 3.46364, so $\mbox{peri}(\mathcal{C}(X)) \leq 3.46364$. Therefore,
\[ \left| \mu(X') - \mu(X) \right| \leq 2.44916 \: d_1 + 0.49993 \: d_2 \]
by ignoring second order terms. \qed
\end{proof}

\begin{figure}[htbp]
	\centering
	\subfigure[3D surface plot]{\label{t01side} \includegraphics[width=4in]{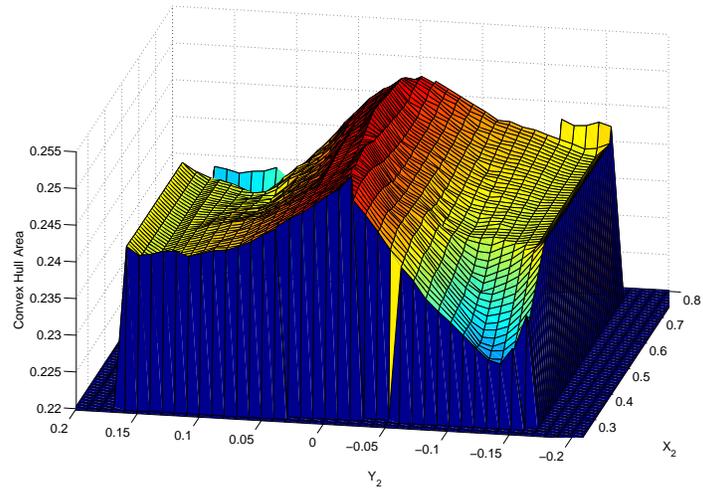}} $\ $
	\subfigure[Top view]{\label{t01top} \includegraphics[width=4in]{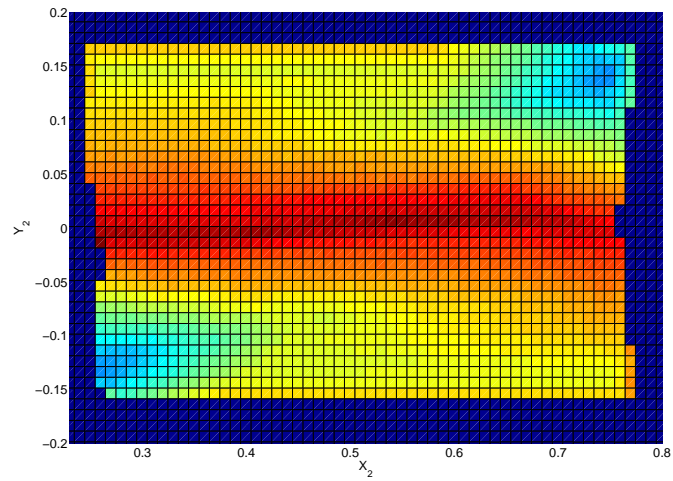}} 
	\caption{Surface plot for $d_1 = d_2 = 0.01$.}
	\label{diagram}
\end{figure}

Now we present the grid-search's result. In the first run, we set both grid sizes $d_1$ and $d_2$ to be $0.01$ and kept track of the optimal possible configuration for each value of $(x_2, y_2)$, the centroid of the triangle $\mathcal{T}$. The surface plot in Figure~\ref{t01side} shows the area of the optimal configuration as a function of $(x_2, y_2)$.

As a quick check, note that the surface has $180^\circ$ rotational symmetry around $(0.5, 0)$ which is what we expected because placing the triangle $\mathcal{T}$ at $(x_2, y_2)$ is the same as placing it at $(1-x_2, -y_2)$. This feature is clearer seen in Figure~\ref{t01top} which is the top-view of our surface plot in Figure~\ref{t01side}. Moreover, the optimal configurations seem to be clustered nicely in one place.

Based on the plot, we heuristically focused our grid-search algorithm to the region bounded by $(x_2, y_2) \in [0.7, 0.77] \times [0.1, 0.17]$ with finer grid size $d_1 = 0.001$. In parallel, we also ran the grid-search algorithm to find the optimal possible configuration for each value of $(x_1, y_1)$ -- the centroid of the square $\mathcal{S}$ -- in order to simultaneously reduce the scope of $(x_1, y_1)$ we are considering.

\begin{figure}[htbp]
  \begin{center}
    \includegraphics[width=2.5in]{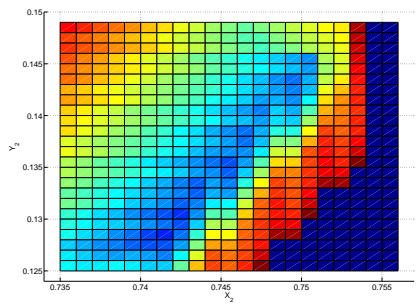}
    \caption{Top-view of the 3D surface plot for $d_1 = 0.001, d_2 = 0.01$.}\label{t001top}
  \end{center}
\end{figure}

\begin{figure}[htbp]
  \begin{center}
    \includegraphics[width=2.5in]{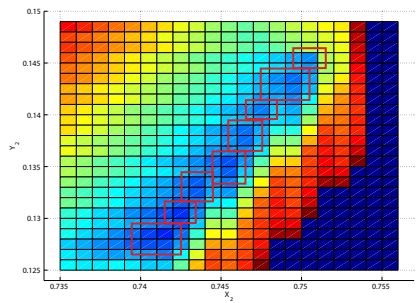}
    \caption{Heuristic scheme for zooming in on the plot in Figure~\ref{t001top}. Each red box represents a region in which we ran independent search.}\label{t001div}
  \end{center}
\end{figure}

The surface plot in Figure~\ref{t001top} reveals that there are indeed multiple grids contributing to small configurations and therefore we zoomed in on each of them with finer grid size $d_2 = 0.0001$ according to the scheme in Figure~\ref{t001div}.

As a result we found an optimal configuration with area $0.227628$ and parameters $(0.6625, 0.1895, 1.30829, 0.7415, 0.1305, 1.63299)$. The smallest step size we were able to run the algorithm is $d_1 = 0.001$ and $d_2 = 0.0001$ . The approximation error of our grid-search is at most $0.0025$ according to Proposition~\ref{stepbound}. Moreover, when we drew small configurations from each of the grid-search runs according to the scheme in Figure~\ref{t001div}, we observed that those configurations are close to a configuration with special features. We conjectured that these are features of a minimal configuration as stated below.

%\begin{figure}[htbp]
%\begin{center}
%\begin{pspicture}(-3,-.3)(3,2.7) \psset{xunit=6cm,yunit=6cm, runit=6cm}
%  \psline(-0.5,0)(.5,0)
%  \pspolygon(0.223556,0.418617)(0.000955,-0.029)(.499989,0.00198144)
% \pspolygon(.223663,.417128)(-0.065128,.250663)(0.101337,-0.0381282)(.390128,.128337)
%\end{pspicture} 
%\end{center}
%  \caption{A configuration with area of 0.227628 obtained from the algorithm.}\label{champ2}
% \end{figure}

\begin{conj} Considering only the case $y_1 \geq 0$, the optimal configuration of $\mathcal{L}$, $\mathcal{T}$ and $\mathcal{T}$ with minimal convex hull area must have the following properties:

(1) The right-most vertex of $\mathcal{T}$ coincides with the point $(1, 0)$ of $\mathcal{L}$

(2) The top-most vertices of $\mathcal{S}$ and $\mathcal{T}$ coincide.
\end{conj}

Our final result is a special grid-search on all configurations satisfying the conjecture. The only parameters are the angle of $\mathcal{T}$ pivoting around $(1, 0)$ and the angle of $\mathcal{S}$ pivoting around the top-most vertex of $\mathcal{T}$. Setting the step size for both parameters to $0.0000001$, we discovered an optimal configuration with area of $0.22758966937711944$ (see Figure~\ref{champ}) which is $0.0001$ more than our improved lower bound. Its parameters are (0.6605, 0.1878, 1.3077, 0.741, 0.1274, 1.6373).

\end{document}